\documentclass[12pt]{article}
\usepackage{amsmath}
\usepackage{amssymb}
\usepackage{amscd}
\usepackage{amsthm}
\theoremstyle{plain}
\newtheorem{Thm}{Theorem}[section]
\newtheorem{Conj}[Thm]{Conjecture}
\newtheorem{Cor}[Thm]{Corollary}
\newtheorem{Lem}[Thm]{Lemma}
\newtheorem{Prob}[Thm]{Problem}
\theoremstyle{definition}
\newtheorem{defn}[Thm]{Definition}
\theoremstyle{remark}

\numberwithin{equation}{section}
\def\Rat{\mathop{\rm Rat}}
\def\Exc{\mathop{\rm Exc}}
\title{\textbf{Tsuji's Numerically Trivial Fibrations and Abundance}}
\author{ \textbf{Shigetaka Fukuda}}
\date{\empty}

\begin{document}
\maketitle
\begin{abstract}
We apply Tsuji's theory of numerically trivial fibrations to the abundance 
problem.
\end{abstract}

In this article every algebraic variety is defined over the field of complex 
numbers 
$\mathbf{C}$ and the standard notation and terminology of the minimal model 
theory (see \cite{Utah}) is used.

We note the famous abundance conjecture:
\begin{Conj}[Abundance]\label{ab}
Let $X$ be an $n$-dimensional projective variety with only 
$\mathbf{Q}$-factorial terminal singularities such 
that the canonical divisor $K_X$ is nef.
Then $K_X$ is semi-ample.
\end{Conj}

The conjecture above is based on the classification philosophy that a minimal 
algebraic variety should be decomposed into varieties of parabolic type (whose 
canonical divisors are numerically trivial) and a variety of hyperbolic type 
(whose canonical divisor is ample).
The log abundance conjecture has been considered by experts to make possible the 
inductive treatment for the abundance conjecture on dimension:
\begin{Conj}[Log Abundance for klt pairs]\label{logab}
Let $(X,\Delta)$ be an $n$- dimensional projective variety with only 
Kawamata log terminal (klt) singularities such 
that the log canonical divisor $K_X + \Delta$ is nef.
Then $K_X + \Delta$ is semi-ample.
\end{Conj}

To state the main theorem, we make two definitions:
\begin{defn}
We say that a property $P(v)$ holds for a {\em very general} point $v$ of an 
algebraic variety $V$ if there exist countably many closed proper subsets 
$G_i$ of $V$ such that the property $P(t)$ holds for every closed point $t \in V 
\setminus (\bigcup G_i)$.
\end{defn}

\begin{defn}[Tsuji \cite{Tsuji}]
For a nef $\mathbf{Q}$-divisor $D$ on an algebraic variety $X$ and a point $x \in 
X$, we set
$N(D;x):= \max \{\dim V ; V$ is a subvariety of $X$ such that $x \in V$ and that 
$D \vert_V$ is numerically trivial$\}$
and
$N(D,X) := \max \{ d \in \mathbf{N} ; N(D;x) \geq d$ for a very general point $x 
\in X \}$.
\end{defn}

The following three problems are well known to experts, concerning the abundance 
conjecture:
\begin{Prob}\label{lmmc}
The log minimal model conjecture for $\mathbf{Q}$-factorial Kawamata log terminal 
(klt) pairs.
\end{Prob}
\begin{Prob}\label{ramif}
Assume that $(X,\Delta)$ is $\mathbf{Q}$-factorial and klt, that $f:X \to S$ is a 
morphism between normal projective varieties with only connected fibers and that 
$K_X+\Delta \sim_{\mathbf{Q}} f^* A$ for some $\mathbf{Q}$-Cartier 
$\mathbf{Q}$-divisor $A$.
Is $K_X+\Delta$ $\mathbf{Q}$-linearly equivalent to $f^* (K_S + \Delta_S)$ for 
some klt pair $(S,\Delta_S)$?
\end{Prob}

Problem \ref{ramif} was treated in Fujino \cite{Fujino}.

\begin{Prob}\label{nnb}
Assume that $(X,\Delta)$ is $\mathbf{Q}$-factorial, projective and klt and that 
$K_X+\Delta$ is nef and $N(K_X+\Delta,X)=0$.
Is $K_X+\Delta$ semi-ample?
(i.e.\ When $K_X+\Delta$ is {\bf numerically ``nef and big"}, is it nef and big?)  
\end{Prob}

Now we state our main theorem:
\begin{Thm}\label{MainT}
If all the answers to Problems \ref{lmmc}, \ref{ramif} and \ref{nnb} are 
affirmative, then Conjecture \ref{logab} is true and so is the abundance 
conjecture \ref{ab}.
\end{Thm}

The main tool for the proof of Theorem \ref{MainT} is Tsuji's theory of 
numerically trivial fibrations:
\begin{Thm}[Tsuji \cite{Tsuji} (cf.\ \cite{BCEKPRSW})]\label{ntf}
Let $L$ be a nef $\mathbf{Q}$-Cartier divisor on a normal projective variety $X$.
Then there exist projective varieties $X_1$ and $Y$ and morphisms $\mu : X_1 \to 
X$ and $f: X_1 \to Y$ with 
the following properties:

(1) $\mu$ is birational.

(2) $f$ is surjective and the function field $\Rat Y$ is algebraically closed in 
$\Rat X_1$.

(3) $\dim Y = \dim X - N(L, X)$.

(4) $X_1 \subset X \times Y$, the morphisms $\mu$ and $f$ are natural projections 
and $\mu^* L \vert_{f^{-1}(y)}$ is numerically trivial for a very general point 
$y \in Y$.

(5) $\mu^{-1}(\mu (f^{-1}(y_0))) = f^{-1}(y_0)$ for some closed point $y_0 \in 
Y$.
\end{Thm}

In Section \ref{Prel}, we collect lemmas and notions.
In Section \ref{PMT}, we prove Theorem \ref{MainT}.
In Section \ref{Appl}, we state some applications.

\section{Preliminaries}\label{Prel}

It is known that the (logarithmic) Kodaira dimension is invariant under the (log) 
minimal model procedure:

\begin{Lem}[\cite{KMM}]\label{discrep}
(1) Let $(Y_{lm},\Delta_{lm})$ be a (relative) log minimal model for a 
$\mathbf{Q}$-factorial klt pair $(Y, \Delta)$.
We consider common resolutions $\gamma_1$ and $\gamma_2$ that are projective 
morphisms:
$$
\begin{CD}
(Y, \Delta) @<\text{$\gamma_1$}<< W @>\text{$\gamma_2$}>> (Y_{lm},\Delta_{lm})
\end{CD}
$$
Then $K_W + {\gamma_1}^{-1}_* \Delta + (1-\epsilon) E \geq {\gamma_1}^* (K_Y + 
\Delta) \geq {\gamma_2}^* (K_{Y_{lm}} + \Delta_{lm})$ for some sufficiently small 
positive rational number $\epsilon$, where $E$ is the reduced divisor composed of 
$\gamma_1$-exceptional prime divisors.

(2) Let $M$ be a (relative) minimal model for a variety $X_2$ with only 
$\mathbf{Q}$-factorial terminal singularities.
We consider common resolutions $\gamma_1$ and $\gamma_2$ that are projective 
morphisms:
$$
\begin{CD}
X_2 @<\text{$\gamma_1$}<< W @>\text{$\gamma_2$}>> M
\end{CD}
$$
Then $K_W \geq {\gamma_1}^* K_{X_2} \geq {\gamma_2}^* K_M$.
\end{Lem}

In this paper we denote the Iitaka dimension by $\kappa$ and the numerical Iitaka 
dimension for a nef $\mathbf{Q}$-divisor by $\nu$.
Nakayama generalized the concept of numerical Iitaka dimension for a nef 
$\mathbf{Q}$-divisor to the case of $\mathbf{Q}$-divisors which are not 
necessarily nef:

\begin{defn}[Nakayama \cite{Ny97}]
For a $\mathbf{Q}$-divisor $D$ and an ample divisor $A$ on a nonsingular 
projective variety $X$, we define $\kappa_{\sigma} (D,A):= \max \{ k \in 
(\mathbf{Z} \cup \{ - \infty \}); \text{there exists a positive number } c \text{ 
such that } \dim H^0 (X, {\cal O}_X (mD+A))$ $\geq c m^k \text{ for every 
sufficiently large and divisible integer } m \}$ and Nakayama's 
$\sigma$-dimension $\kappa_{\sigma} (D,X):= \max \{ \kappa_{\sigma} (D,A) ; A 
\text{ is an ample divisor on } X \}$.
\end{defn}

\begin{Lem}[Nakayama \cite{Ny97}]\label{KS}
Let $D$ be a $\mathbf{Q}$-Cartier $\mathbf{Q}$-divisor on a normal projective 
variety 
$X$, $f:Z \to X$ a surjective morphism from a nonsingular projective variety and 
$E$ an $f$-exceptional effective $\mathbf{Q}$-divisor.

(1) In the case where $X$ is nonsingular, $\kappa_{\sigma} (D,X) \geq 0$ if and 
only if $D$ is pseudo-effective.

(2) In the case where $X$ is nonsingular, $\kappa_{\sigma} (f^* D,Z) = 
\kappa_{\sigma} (D,X)$.

(3) $\kappa_{\sigma} (f^* D + E, Z) = \kappa_{\sigma} (f^* D,Z)$ if $f$ is 
birational.

(4) In the case where  $X$ is nonsingular and $D$ is nef, $\kappa_{\sigma} (D, X) 
= \nu (D, X)$.
\end{Lem}

\section{Proof of Theorem \ref{MainT}}\label{PMT}
We prove the theorem by induction on dimension.
We may assume that $X$ is $\mathbf{Q}$-factorial by taking a $K_X + 
\Delta$-crepant 
$\mathbf{Q}$-factorialization due to Problem \ref{lmmc}.

{\it Case 1:} $N(K_X + \Delta, X) = n$ (i.e.\ $K_X + \Delta$ is numerically 
trivial).
We use Kawamata's argument (\cite{Ka85b}, the proof of Theorem 8.2).

{\it Subcase 1(a):} $q(X)=0$.
In this subcase, the numerical triviality implies the $\mathbf{Q}$-linear 
triviality.

{\it Subcase 1(b):} $q(X) \geq 1$.
We have the Albanese morphism $f:X \to S$.
Consider the Stein factorization:
$$
\begin{CD}
f:X @>\text{$g$}>> X' @>\text{$h$}>> f(X)
\end{CD}
$$
If $\dim f(X) < \dim X$ then, from the induction hypothesis $(K_X + \Delta) 
\vert_{g^{-1} (x')} = K_{g^{-1} (x')} + \Delta \vert_{g^{-1} (x')}$ is 
$\mathbf{Q}$-linearly trivial for a general point $x' \in X'$ and thus, from 
Nakayama (\cite{Ny86}, Theorem 5), $K_X + \Delta \sim_{\mathbf{Q}} g^* A$ for 
some $\mathbf{Q}$-Cartier $\mathbf{Q}$-divisor $A$ on $X'$.
Therefore Problem \ref{ramif} and the induction hypothesis imply the assertion.
If $\dim f(X) = \dim X$, then the Kodaira dimension of a smooth model of $f(X)$ 
is non-negative from the theory of Abelian varieties, thus so is the Kodaira 
dimension of a smooth model of $X$.
Consequently we have $\kappa (K_X + \Delta) \geq 0$.

{\it Case 2:} $0 < N(K_X + \Delta, X) < n$.
From Tsuji (Theorem \ref{ntf}), there exist projective 
varieties $X_1$ and $Y$ and morphisms $\mu : X_1 \to X$ and $f: X_1 \to Y$ with 
the following properties:

(1) $\mu$ is birational.

(2) $f$ is surjective and the function field $\Rat Y$ is algebraically closed in 
$\Rat X_1$.

(3) $\dim Y = \dim X - N(K_X + \Delta, X)$.

(4) $X_1 \subset X \times Y$, the morphisms $\mu$ and $f$ are natural projections 
and $\mu^* 
(K_X + \Delta) \vert_{f^{-1}(y)}$ is numerically trivial for a very general point 
$y \in Y$.

(5) $\mu^{-1}(\mu (f^{-1}(y_0))) = f^{-1}(y_0)$ for some closed point $y_0 \in 
Y$.

From (5) and Zariski's Main Theorem, there exists an open subset $X_0$ of $X$ 
such that $\mu^{-1}(X_0) \supset f^{-1}(y_0)$ and  $X_0 \cong (\text{the 
normalization of } \mu^{-1}(X_0))$.
Then $X_0 \cong \mu^{-1}(X_0)$.
Put $V:= Y \setminus f (X_1 \setminus \mu^{-1} (X_0))$ and $U:= \mu (f^{-1} 
(V))$.
We note that $V$ and $U$ are non-empty and open.
Then $U \cong \mu^{-1} (U) = f^{-1} (V)$.
For a very general point $y \in Y$, the pair $(f^{-1}(y), \mu^* \Delta 
\vert_{f^{-1} (y)})
$ is with only klt singularities and $K_{f^{-1}(y)} + \mu^* \Delta \vert_{f^{-1} 
(y)}= \mu^* (K_X + \Delta) \vert_{f^{-1}(y)}$ is numerically trivial and thus, 
from the induction hypothesis, $\mathbf{Q}$-linearly trivial.

This paragraph is a variant of Nakayama (\cite{Ny91}, Appendix A).
Consider the following commutative diagram:
$$
\begin{CD}
X_2 @= X_2 \\
@V\text{$\rho$}VV @VV\text{$g$}V \\
X_1 @>>\text{$f$}> Y \\
@V\text{$\mu$}VV \\
X
\end{CD}
$$
where $\mu \rho$ is a log resolution of $(X,\Delta)$ and $X_2$ is a projective 
variety.
According to Problem \ref{lmmc} we can apply the relative log minimal model 
program for $g: (X_2, (\mu \rho)^{-1}_* \Delta + (1 - \epsilon) E_2) \to Y$ where 
$E_2$ is the reduced divisor composed of $\mu \rho$-exceptional prime divisors 
and 
where $\epsilon$ is a sufficiently small positive rational number.
We end up with the relative log minimal model $h: (M, \Delta_M) \to Y$.
Here $K_M + \Delta_M$ is $h$-nef and $0 = \kappa (K_{f^{-1} (y)} + \mu^* \Delta 
\vert_{f^{-1} (y)}, f^{-1} (y)) = \kappa (K_{g^{-1} (y)} + (\mu \rho)^{-1}_* 
\Delta \vert_{g^{-1} (y)} + (1-\epsilon) E_2 \vert_{g^{-1} (y)}, g^{-1} (y)) = 
\kappa (K_{h^{-1}(y)} + \Delta_M \vert_{h^{-1}(y)}, h^{-1}(y))$ for a very 
general point $y \in Y$ from Lemma \ref{discrep}  (we note $\Exc (\rho 
\vert_{g^{-1}(y)}) = \Exc (\rho) \vert_{g^{-1}(y)}$).
Then $(K_M + \Delta_M) \vert_{h^{-1}(y)} = K_{h^{-1}(y)} + \Delta_M \vert_{h^{-1}
(y)}$ is $\mathbf{Q}$-linearly trivial from the induction hypothesis.
Thus from Nakayama (\cite{Ny86}, Theorem 5) (see also \cite{KMM}, 6-1-11), we 
obtain a projective variety $Y_2 := \mathbf{P} (\bigoplus_{l \geq 0} h_* {\cal O}
_M (l (K_M + \Delta_M))$, a natural morphism $p: M \to Y_2$ and a 
$\mathbf{Q}$-Cartier $\mathbf{Q}$-divisor $A$ on $Y_2$ such that $K_M + \Delta_M 
\sim_{\mathbf{Q}} p^* A$.
From Problem \ref{ramif}, $K_M + \Delta_M \sim_{\mathbf{Q}} p^* (K_{Y_2} + 
\Delta_2)$, for some klt pair $(Y_2,\Delta_2)$.

We consider the following commutative diagram:
$$
\begin{CD}
X_2 @<\text{$\alpha$}<< W_1 @= W_1 @>\text{$\delta$}>> W_2 
@>\text{$\gamma_2$}>> Y_{lm} \\
@V\text{$g$}VV @V\text{$\beta$}VV @VVV @VV\text{$\gamma_1$}V \\
Y @<<\text{$h$}< M @>>\text{$p$}> Y_2 @<<\text{q}< Y_Q
\end{CD}
$$
where the morphisms $\alpha$ and $\beta$ are projective, they are common 
resolutions of 
$X_2$ and $M$, the morphism $q:(Y_Q,\Delta_Q) \to (Y_2,\Delta_2)$ is a 
$K_X+\Delta$-crepant $\mathbf{Q}$-factorialization of $(Y_2,\Delta_2)$ according 
to Problem 
\ref{lmmc} (i.e.\ $(Y_Q,\Delta_Q)$ is projective, $\mathbf{Q}$-factorial and klt 
and $q^* (K_{Y_2} + \Delta_2) = K_{Y_Q} + \Delta_Q$), the log variety $(Y_{lm}, 
\Delta_{lm})$ is a log minimal model for $(Y_Q,\Delta_Q)$ according to Problem 
\ref{lmmc}, the morphisms $\gamma_1$ and $\gamma_2$ are projective and they are 
common 
resolutions of $Y_Q$ and $Y_{lm}$.

From Lemmas \ref{discrep} and \ref{KS} we derive the equalities $\nu (K_X + 
\Delta,X) = \nu ( \rho^* \mu^* (K_X + \Delta),X_2) = \kappa_{\sigma} ( \rho^* 
\mu^* (K_X + \Delta),X_2) = \kappa_{\sigma} (K_{X_2} + (\mu \rho)^{-1}_* \Delta + 
(1-\epsilon) E_2, X_2) = \kappa_{\sigma} (\alpha^* (K_{X_2} + (\mu \rho)^{-1}_* 
\Delta + (1-\epsilon) E_2), W_1) = \kappa_{\sigma} (K_{W_1}+ \alpha^{-1}_* ((\mu 
\rho)^{-1}_* \Delta) + \alpha^{-1}_* ((1-\epsilon) E_2) + E_W, W_1)= 
\kappa_{\sigma} (\beta^* (K_M + \Delta_M), W_1) = \kappa_{\sigma} (\beta^* p^* 
(K_{Y_2} + \Delta_2), W_1) = \kappa_{\sigma} (\delta^* {\gamma_1}^* q^* (K_{Y_2} 
+ \Delta_2), W_1) = \kappa_{\sigma} ({\gamma_1}^*(K_{Y_Q} + \Delta_Q), W_2)$ $= 
\kappa_{\sigma} (K_{W_2} + {\gamma_1}^{-1}_* \Delta_Q + F_W, W_2) = 
\kappa_{\sigma} ({\gamma_2}^* (K_{Y_{lm}} + \Delta_{lm}), W_2) = \nu 
({\gamma_2}^* (K_{Y_{lm}} + \Delta_{lm}), W_2) = \nu (K_{Y_{lm}} + \Delta_{lm}, 
Y_{lm})$, where $E_W$ and $F_W$ are the reduced divisors composed of 
$\alpha$-exceptional and $\gamma_1$-exceptional prime divisors respectively.
Because $\kappa (K_{Y_{lm}} + \Delta_{lm}, Y_{lm}) = \nu (K_{Y_{lm}} + 
\Delta_{lm}, Y_{lm})$ from the induction hypothesis, the equalities above implies 
that $\kappa (K_X + \Delta, X) = \nu (K_X + \Delta, X)$.
Thus Kawamata \cite{Ka85a} implies the assertion.

{\it Case 3:} $N(K_X + \Delta, X) = 0$.
This is Problem \ref{nnb}.
\qed

\section{Applications}\label{Appl}

\begin{Thm}
Let $X$ be an $n$-dimensional projective variety with only terminal singularities 
such that the canonical divisor $K_X$ is nef.
Assume that the minimal model conjecture is true in dimension $n$ and that the 
log minimal model and the log abundance conjectures for $\mathbf{Q}$-factorial 
klt pairs are true in dimension $n-1$.
If $N(K_X,X)=1$, then the canonical divisor $K_X$ is semi-ample.
\end{Thm}

\begin{proof}
We note that the number $N(K_X, X)$ is invariant under the procedure of 
$K_X$-crepant $\mathbf{Q}$-factorialization.
Problem \ref{ramif} is solved affirmative by Nakayama (\cite{Ny87}, 0.4) in the 
case where $X$ is with only terminal singularities, $\Delta = 0$ and $\dim X - 
\dim S = 1$.
Thus the proof of Theorem \ref{MainT} works.
\end{proof}

If the four-dimensional minimal model strategy of Shokurov is completed, the 
following becomes important:

\begin{Cor}
Let $X$ be a projective fourfold with only terminal singularities such that the 
canonical divisor $K_X$ is nef.
Assume that the minimal model conjecture is true in dimension four.
If $N(K_X,X)=1$, then $K_X$ is semi-ample.
\end{Cor}

\begin{proof}
The log minimal model and the log abundance conjectures for 
$\mathbf{Q}$-factorial klt pairs in dimension three are proved by Shokurov 
\cite{Sh} and by Keel, Matsuki and McKernan \cite{KeMaMc} respectively.
\end{proof}

\begin{Thm}
Assume that $(X,\Delta)$ is projective, $\mathbf{Q}$-factorial and klt and that 
$K_X+\Delta$ is numerically trivial.
If the answer to Problem \ref{ramif} is affirmative, then $K_X+\Delta$ is 
$\mathbf{Q}$-linearly trivial.
\end{Thm}

\begin{proof}
See Case 1 in the proof of Theorem \ref{MainT}.
\end{proof}

The following is a variant of a result of Bauer, Campana, Eckl, Kebekus, 
Peternell, Rams, Szemberg and Wotslaw:

\begin{Thm}[cf.\ \cite{BCEKPRSW}, 2.4.4]
Let $(X,\Delta)$ be a projective klt $n$-fold such that $K_X + \Delta$ is nef and 
$N(K_X + \Delta,X) = n-1$. 
Assume that every numerically trivial log canonical divisor on a projective klt 
log variety is $\mathbf{Q}$-linearly trivial in dimension $n-1$.
Then $K_X + \Delta$ is semi-ample.
\end{Thm}

\begin{proof}
We argue along the lines in the proof of Theorem \ref{MainT}.
From Tsuji (Theorem \ref{ntf}), there exist projective 
varieties $X_1$ and $Y$ and morphisms $\mu : X_1 \to X$ and $f: X_1 \to Y$ with 
the following properties:

(1) $\mu$ is birational.

(2) $f$ is surjective and the function field $\Rat Y$ is algebraically closed in 
$\Rat X_1$.

(3) $\dim Y = 1$.

(4) $X_1 \subset X \times Y$, the morphisms $\mu$ and $f$ are natural projections 
and $\mu^* 
(K_X + \Delta) \vert_{f^{-1}(y)}$ is numerically trivial for a very general point 
$y \in Y$.

(5) $\mu^{-1}(\mu (f^{-1}(y_0))) = f^{-1}(y_0)$ for some closed point $y_0 \in 
Y$.

From (5), the rational map $f \cdot \mu^{-1}: X \cdots \to Y$ is a morphism 
(\cite{BCEKPRSW}, 2.4.4).
We consider the Stein factorization:
$$
\begin{CD} 
f \cdot \mu^{-1} :X @>\text{$p$}>> Y_2 @>>> Y.
\end{CD}
$$
From the assumption, $(K_X + \Delta) \vert_{p^{-1} (y_2)}$ is 
$\mathbf{Q}$-linearly trivial for a very general point $y_2 \in Y_2$.
Thus $K_X + \Delta \sim_{\mathbf{Q}} p^* A$ for some $\mathbf{Q}$-Cartier 
$\mathbf{Q}$-divisor $A$ on $Y_2$ by Nakayama (\cite{Ny86}, Theorem 5).
From the fact that $\nu (K_X + \Delta, X) > 0$, we have $\deg A > 0$.
Thus $K_X + \Delta$ is semi-ample.
\end{proof}

\begin{Cor}[cf. \cite{BCEKPRSW}, 2.4.4]
Let $(X,\Delta)$ be a projective klt fourfold such that $K_X + \Delta$ is nef and 
that $N(K_X + \Delta,X) = 3$.
Then the log canonical divisor $K_X + \Delta$ is semi-ample.
\end{Cor}

\begin{proof}
The assumption of the theorem is satisfied from Keel, Matsuki and McKernan 
\cite{KeMaMc}.
\end{proof}

\bigskip
Faculty of Education, Gifu Shotoku Gakuen University

Yanaizu-cho, Gifu 501-6194, Japan

fukuda@ha.shotoku.ac.jp

\end{document}